\documentclass[12pt,leqno]{article}
\usepackage{graphicx}
\newcommand{\PP}{{\rm I\!P}}

\newcommand{\R}{{\rm I\!R}}

\newcommand{\EE}{{\rm I\kern-2pt E}}
\newcommand{\RR}{{\rm I\kern-2pt R}}
\newcommand{\sumi}{\sum_{i=1}^n}

\newcommand{\Frac}{\frac{1}{2}}
\newcommand{\lam}{\lambda}
\newcommand{\sig}{\sigma}
\newcommand{\Ome}{\Omega}
\newcommand{\Fc}{{\cal F}}
\newcommand{\del}{\delta}
\newcommand{\Del}{\Delta}
\newcommand{\ome}{\omega}
\newcommand{\varp}{\varphi}
\newcommand{\Sig}{\Sigma}
\newtheorem{thm}{Theorem}

\newtheorem{pro}[thm]{Proposition}

\newtheorem{dfn}[thm]{Definition}

\newcommand{\alp}{\alpha}
\newcommand{\bet}{\beta}
\newcommand{\Gam}{\Gamma}
\title{Viscoelasticity and L\'evy processes\\ Visco-\'elasticit\'{e} et Processus de L\'{e}vy}
\author{Nicolas BOULEAU}\date{\it Ecole des Ponts, ParisTech.}
\begin{document}

\maketitle

\noindent{\bf Abstract} We show that the linear viscoelastic materials, and more generally the physical phenomena
 to which Biot's relaxation theory is relevant, can be put in correspondance with
 the laws of processes with independent increments.

In the one dimensional case this correspondence is one to one with subordinators
 and gives rise naturally to a conjugation relation on subordinators.\\



{\bf PRESENTATION}\\

Certaines \'equations de la physique et les probl\`emes aux limites associ\'es
admettent une interpr\'etation probabiliste. C'est le cas de l'\'equation de la
chaleur associ\'ee au mouvement brownien, ou des \'equations elliptiques
lin\'eaires du second ordre associ\'ees aux processus de diffusion, ou encore
des op\'erateurs lin\'eaires int\'egro-diff\'erentiels v\'erifiant le principe
du maximum positif auxquels sont associ\'es des processus de Markov avec sauts.
Ceci est \'etudi\'e par la th\'eorie probabiliste du potentiel (cf par exemple
[4], [11], [12]). Le cas discontinu est une extension math\'ematique du cas
continu et correspond \`a des probl\`emes non locaux plus rarement rencontr\'es
par l'ing\'enieur. De nombreuses autres \'equations (cf [14]) poss\`edent aussi maintenant 
des
interpr\'etations probabilistes.

Nous montrons ici que les ph\'enom\`enes visco\'elastiques lin\'eaires sont
susceptibles d'une interpr\'etation probabiliste par des processus discontinus.
Cette interpr\'etation est analogue aux pr\'ec\'edentes, en ce sens, que les
grandeurs physiques y apparaissent comme l'esp\'erance de fonctionnelles du
processus. N\'eanmoins, ici le temps qui
r\'egit le mouvement de la particule probabiliste n'est pas le temps qui
r\'egit le ph\'enom\`ene physique (comme c'est le cas lorsqu'on associe au
 semigroupe de la chaleur un mouvement brownien). Par ailleurs, cette correspondance peut
\^etre faite de deux fa\c{c}ons qui sont duales et ceci conduit \`a d\'egager
la notion de couple conjugu\'e. La pr\'{e}sente r\'{e}daction est une version d\'{e}taill\'{e}e de [6].

\noindent Une fonction $\phi : \RR_+ \rightarrow \RR_+$ est une
fonction de Bernstein si $\phi$ est $C^\infty$,  $\phi \geq 0$ et $ (-1)^n D^n \phi \leq
0\;\;\;\forall n\geq 1$. Ces fonctions interviennent
 dans la caract\'{e}risation des
semi-groupes de convolution de probabilit\'{e}s sur $\RR_+$ avec
lesquels elles sont en correspondance biunivoque (cf [1])
ainsi donc qu'avec les processus \`{a} accroissements
ind\'{e}pendants croissants ou subordinateurs. Ces fonctions
forment un c\^{o}ne stable par composition dont la structure est
assez complexe (cf [7], [8]). Elles interviennent aussi, ceci est
corr\'{e}latif, dans la th\'{e}orie du calcul symbolique sur les
g\'{e}n\'{e}rateurs infinit\'{e}simaux de semi-groupes de Markov
(cf [1], [15], [18]). 

Il est  remarquable que ce soit exactement cette classe de fonctions
qui intervienne dans les ph\'{e}nom\`{e}nes de relaxation tels
que la visco\'{e}lasticit\'{e} lin\'{e}aire qui sont r\'{e}gis
par une math\'{e}matique {\it a priori} sans rapport avec les
objets pr\'{e}c\'{e}dents mais qui se trouve mise en connexion
avec eux par le raisonnement par variables cach\'{e}es. Plus
pr\'{e}cis\'{e}ment: Les ph\'{e}nom\`{e}nes de relaxation tels
que la visco\'{e}lasticit\'{e} lin\'{e}aire sont d\'{e}crits par
des op\'{e}rateurs dissipatifs et sont donc li\'{e}s d'apr\`{e}s
le th\'{e}or\`{e}me de Lumer-Phillips ([23] p.250) \`{a}
des semi-groupes \`{a} contraction mais ces semi-groupes n'ont
aucune interpr\'{e}tation probabiliste en g\'{e}n\'{e}ral car
ils n'op\`{e}rent pas positivement sur les fonctions; dans cette
situation c'est le raisonnement par variables cach\'{e}es qui
fait appara\^{\i}tre des objets li\'{e}s \`{a} des semi-groupes
de Markov.

\section{Thermodynamique et visco\'{e}lasticit\'{e}}
\subsection{Notations}

Les mod\`{e}les visco-\'{e}lastiques les plus simples sont les  mod\`{e}les
de Maxwell (pour un liquide visco-\'{e}lastique):

\begin{center}
\includegraphics[width=3.4in]{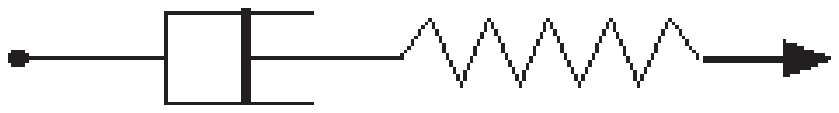}
\end{center}

\noindent et de Kelvin-Voigt (pour un solide visco-\'{e}lastique):

\begin{center}
\includegraphics[width=3in]{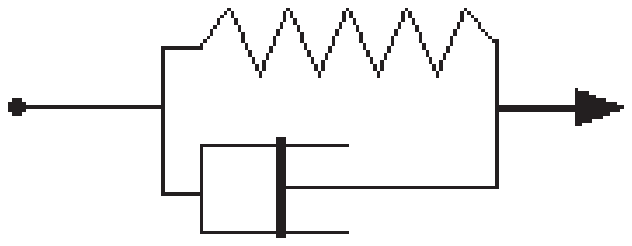}
\end{center}

Lorsqu'un tel syst\`{e}me est soumis \`{a} une force d\'{e}pendant du temps $Q(t)$,
 la d\'{e}formation $q(t)$ qui en r\'{e}sulte d\'{e}pend lin\'{e}airement de $Q$ et si les propri\'{e}t\'{e}s
 du corps ne varient pas avec le temps cette transformation est connue si on conna\^{\i}t la r\'{e}ponse 
pour l'impulsion unit\'{e} appel\'{e}e  {\it r\'{e}ponse
impulsionnelle} $f(t)$ du syst\`{e}me dont la r\'{e}ponse g\'{e}n\'{e}rale est alors 
$$q(t)=\int_{[0,t]}f(t-\tau) \;dQ(\tau).$$ Pour un mod\`{e}le de Maxwell
on a $$\dot{q}=\frac{\dot{Q}}{G}+\frac{Q}{\eta}$$ o\`{u} $G$ est la raideur du ressort et $\eta$
la viscosit\'{e} de l'amortisseur, et donc $$f(t)=(\frac{1}{G}+\frac{t}{\eta})1_{\{t\geq 0\}}.$$
Par ailleurs on peut au contraire imposer une d\'{e}formation au syst\`{e}me et mesurer la force qui assure l'\'{e}quilibre. 
Pour une d\'{e}formation qui passe de 0 \`{a} 1 \`{a} l'instant 0 puis ne varie plus, la force mesur\'{e}e est la {\it 
fonction de relaxation} $r(t)$ du syst\`{e}me dont la r\'{e}ponse g\'{e}n\'{e}rale est alors 
$$Q(t)=\int_{[0,1]} r(t-\tau)\;dq(\tau)$$ pour le mod\`{e}le de Maxwell $r(t)=Ge^{-\frac{G}{\eta}t}1_{\{t\geq 0\}}.$

\subsection{Limitations impos\'{e}es par la thermodynamique}
Les restrictions que la thermodynamique des ph\'enom\`enes irr\'eversibles
impose aux ph\'enom\`enes visco\'elastiques ont \'et\'e \'etudi\'es et
d\'egag\'es par plusieurs auteurs parfois ind\'ependamment (cf [13] [21] [2] [9]
[16] [20]). Nous nous pla\c{c}ons sous les hypoth\`eses de la th\'eorie de Biot
[2] et de m\^eme  que celle-ci les consid\'erations qui suivent peuvent
s'appliquer \`a d'autres ph\'enom\`enes physiques (\'electriques ou chimiques)
pourvu que les hypoth\`eses (Principe d'Onsager, existence de variables
normales, lin\'earit\'e) soient acceptables au moins en premi\`ere
approximation. Pour la commodit\'e, nous emploierons le langage de la
m\'ecanique des solides (cf [20]).

Soit un syst\`eme soumis \`a des forces g\'en\'eralis\'ees $Q_i, i=1,\cdots,n$
et d\'ecrit par des param\`etres g\'eom\'etriques associ\'es $q_i, i=1,\cdots,n$
tels que le travail des forces ext\'erieures, s'\'ecrive $\sumi Q_i dq_i$.

Au voisinage d'une position d'\'equilibre stable o\`u l'on a pris $q_i = 0
\;\;\forall i$, le potentiel thermodynamique du syst\`eme s'\'ecrit~:

\[
     W = \Frac \sum_{i,j} a_{ij} q_i q_j
\]
o\`u la matrice $(a_{ij})$ est sym\'etrique semi-d\'efinie positive (le terme
stable est pris au sens large). En  calculant la variation d'entropie durant un
court intervalle de temps (cf [16] chapitre 13) on obtient, pour des vitesses
$q_i$ suppos\'ees petites et en n\'egligeant les forces d'inertie, que la
puissance dissip\'ee peut s'\'ecrire~:
\[
	    D = \Frac \sum_{i,j} b_{ij} \dot{q}_i \dot{q}_j
\]
o\`u la matrice $(b_{ij})$ est sym\'etrique d'apr\`es le principe d'Onsager et
semi-d\'efinie positive par le second principe de la thermodynamique.
L'\'equation d'\'evolution est alors:

\[
     \frac{\partial D}{\partial \dot{q}_i} + \frac{\partial W}{\partial q_i} =
Q_i \]
soit
\begin{equation}
\label{equa.1}
     Q_i = \sum_j a_{ij} q_j + b_{ij} \dot{q}_j
\end{equation}

Le cas de la visco\'elasticit\'e lin\'eaire est celui o\`u les coefficients
$(a_{ij})$ et $(b_{ij})$ sont constants. Dans ce cas, la relation lin\'eaire
entre les histoires des forces $(Q_i(t))$ et celles des param\`etres $(q_j(t))$
commute avec les translations du temps et est connue par la r\'eponse
impulsionnelle $f_{ij}(t)$ du param\`etre $q_i$ \`a la force $Q_j$.

 Ecrivons l'\'{e}quation (\ref{equa.1}) sous la forme

\begin{equation}
\label{equa.2}
     Aq + B\dot{q} = Q
\end{equation}
o\`{u} $A$ et $B$ sont des matrices $n \times n$ sym\'{e}triques semi-d\'{e}finies
positives
et supposons d'abord que  $B$ est d\'{e}finie positive.
Consid\'{e}rons sur $\R^n$ la structure euclidienne associ\'{e}e \`{a} $B$
 dont le produit scalaire est

\[
     (u,v)_B  =  < u,Bv >  =  {^t}u Bv
\]
( $< .,. >$ \'{e}tant le produit scalaire usuel sur $\R^n$). Pour cette nouvelle
structure euclidienne
l'op\'{e}rateur $B^{-1}A$ est auto-adjoint:

\[
     (u,B^{-1}Av)_B = < u,Av > = < Au,v > = (v,B^{-1}Au)_B
\]
Donc il existe une base $B$-orthonormale  $(\psi_k)_{k=1,\ldots,n}$ sur laquelle
$B^{-1}A$ est diagonale, c'est \`{a} dire 

\[
\begin{array}{lclcl}
     (\psi_j,B^{-1}A \psi_k)_B & = & < \psi_j,A \psi_k > & = &
     \left\{ \begin{array}{ll}
        0 & \mbox{ if } j \neq k\\
        \lam_k \geq 0 & \mbox{ if } j = k
     \end{array}\right.\\
	    & & & &\\
     (\psi_j,\psi_k)_B & = & < \psi_j,B \psi_k > & = &
     \left\{ \begin{array}{ll}
        0 & \mbox{ if } j \neq k\\
        1 & \mbox{ if } j = k
     \end{array}\right.
\end{array}
\]
en particulier

\[
     A \psi_k - \lam_k B \psi_k = 0 \qquad \forall k.
\]
Soit $\hat{q}(\theta)$ $[$resp. $\hat{Q}(\theta)]$ la transform\'{e}e de
Laplace de $q(t)$ $[$resp. $Q(t)]$ $(\hat{q}(\theta) =
\int_0^\infty e^{-\theta t} q(t) dt)$.

L'\'{e}quation (\ref{equa.2}) donne

\begin{equation}
\label{equa.3}
     (A + \theta B) \hat{q} = \hat{Q}.
\end{equation}
Donc, si $\hat{q}(\theta)$ est repr\'{e}sent\'{e}e sur la base $(\psi_k)$,

\begin{equation}
\label{equa.4}
     \hat{q}(\theta) = \sum_{k=1}^{n} \xi_k(\theta) \psi_k
\end{equation}
on a 

\[
     (\lam_k + \theta) \xi_k = ( \hat{Q},\psi_k, )_{B}.
\]

Si nous supposons maintenant que seulement $m$ des $n$ param\`{e}tres $ (m < n)$ sont observables,  i.e. if
$\hat{Q} =
(\hat{Q}_1,\ldots,\hat{Q}_m,0,\ldots,0)$ il vient

\[
     \xi_k = \frac{1}{\lam_k + \theta} \sum_{j=1}^{m}
     \hat{Q}_j(B \psi_{k})_j
\]
et par (4)

\[
     \hat{q}_i = \sum_{j=1}^{m} \hat{Q}_j
     [\sum_{k=1}^{n} \frac{1}{\lam_k + \theta}  (B \psi_{k})_j\psi_{ki}]. 
\]
Ainsi, d\'{e}signant par $f_{ij}(t)$ the r\'{e}ponse de $q_i$ \`{a} l'impulsion unit\'{e} de
$Q_j$, on a 

\[
     \hat{f_{ij}^{'}} = \sum_{k=1}^{n} \frac{1}{\lam_k + \theta} 
     \psi_{k_j} (B \psi_{k})_j\psi_{ki}
\]
et donc

\[
     f_{ij}(t) = \sum_{k=1}^{n} \int_0^t e^{-\lam_ks} ds J_{ij}^{(k)}.
\]

On voit facilement que si des hypoth\`{e}ses plus larges sont prises sur la matrice 
 $B$ autorisant des valeurs propres nulles ou infinies, on obtient la forme plus g\'{e}n\'{e}rale

\begin{equation}
\label{equa.5}
     f_{ij}(t) = \sum_{k=1}^{n} (1 - e^{-\lam_kt}) J_{ij}^{(k)} + t L_{ij} +
     K_{ij}
\end{equation}
o\`{u} les matrices $(J^{(k)}) k = 1,\ldots,n,L,K$ sont sym\'{e}triques  semi-definies 
positives.

Partant de l'\'{e}quation (5) et passant \`{a} la limite on obtient que 
les r\'{e}ponses impulsionnelles des syst\`{e}mes visco-\'{e}lastiques sont

\noindent dans le cas d'un seul param\`{e}tre observable

\begin{equation}
\label{equa.6}
     f(t) = \int_{\RR_+^*} (1 - e^{-\lam t}) \nu(d\lam) + bt + c
\end{equation}
o\`{u} $\nu$ est une mesure $\sig$-finie positive on $\RR_+^* = ]0,\infty[$ telle que
 $\int_0^\infty \frac{x}{1 + x} d\nu(x) < +\infty$, et $b \geq 0$, $c \geq
0$; 

\noindent et dans le cas de  $m$ param\`{e}tres observables de la forme

\begin{equation}
\label{equa.7}
     f_{ij}(t) = \int_{\RR_+^*} (1 - e^{-\lam t}) \nu_{ij}(d\lam) + t L_{ij} +
     K_{ij},
\end{equation}
o\`{u} $\nu = (\nu_{ij})$ est une matrice semi-d\'{e}finie positive de mesures
$\sig$-finies sur $\RR_+^*$ satisfaisant

\[
     \int_{\RR_+^*} \frac{x}{1 + x} d|\nu_{ij}|(x) < +\infty \qquad \forall
     i,j = 1,\ldots,n
\]
et o\`{u} les matrices $L$ and $K$ sont  sym\'{e}triques semi-d\'{e}finies positives.

Ce passage \`{a} la limite peut se faire selon des arguments physiques, en prenant
 l'ensemble des limites simples des fonctions de la forme 
 (5), ou
en consid\'{e}rant un continu visco-\'{e}lastique, les matrices $A$ et $B$ \'{e}tant
 remplac\'{e}es
par des op\'{e}rateurs auto-adjoints non born\'{e}s sur un espace de Hilbert.
 Alors l'argument ci-dessus s'\'{e}tend et donne directement 
les formes (6) et (7) par la repr\'{e}sentation
 spectrale des op\'{e}rateurs auto-adjoints.

R\'eciproquement, les syst\`emes caract\'eris\'es par (7)
v\'erifient les exigences de la thermodynamique et de la stabilit\'e (cf  [20]
annexe XXI-20 ou [16] chapitre 13).\\

\noindent {\bf Remarque}. Dans le cas d'un seul param\`{e}tre observable (6) montre que
 la r\'{e}ponse impulsionnelle g\'{e}n\'{e}rale est une fonction de Bernstein. {\it Il est donc faux de croire (comme il est \'{e}crit dans certains manuels) 
que le groupement en parall\`{e}le d'une suite infinie de mod\`{e}les de Maxwell (ou le groupement en s\'{e}rie 
d'une suite infinie de mod\`{e}les de Kelvin-Voigt) donne le mod\`{e}le rh\'{e}ologique visco-\'{e}lastique 
g\'{e}n\'{e}ral, de m\^{e}me que les mesures discr\`{e}tes ne sont pas les mesures 
g\'{e}n\'{e}rales}. On perd en particulier tous les mod\`{e}les o\`{u} la mesure $\nu$ est absolument continue.
Ainsi, de nombreux mod\`{e}les
sont analytiquement calculables qui ne correspondent pas \`{a} des groupements d'amortisseurs et de ressorts.

\section{Processus \`{a} Accroissements Ind\'{e}pendants Stationnaires}

	Pour la commodit\'e de la suite, tous les processus al\'eatoires seront
index\'es par des lettres grecques.

Nous dirons qu'un processus $(Y_\tau)_{\tau\geq0}$ d\'efini sur un espace de
probabilit\'e $(\Ome,{\cal A},\PP)$ est un P.A.I.S. s'il est continu \`a droite
et tel que $\forall \tau_1 < \tau_2 <\cdots< \tau_n$ les variables $Y_0,
Y_{\tau_1} - Y_0,\cdots,Y_{\tau_n} - Y_{\tau_{n-1}}$ sont ind\'ependantes et si
la loi de $Y_\tau - Y_\sig$ ne d\'epend que de $\tau - \sig$.

On note $\Fc_\tau = \sig(Y_\sig, \sig \leq \tau)$. Un processus $(Y_\tau)$
d\'efini sur $(\Ome,{\cal A},\PP)$, continu \`a droite, \`a valeurs $\RR$, est un
P.A.I.S. si et seulement si (formule de L\'evy-Khintchine) (cf [15] [19]).

\[
     \forall u \in \RR^m, \forall \tau \geq \sig \geq 0 \qquad \qquad \qquad
\]
\begin{equation}
\label{equa.4}
     \EE[exp\{i < u,Y_\tau - Y_\sig >\}|\Fc_\sig]\qquad \qquad \qquad \qquad
\end{equation}
\[
     = exp\{(\tau-\sig)[-\Frac < u,\Sigma u > +i < u,\eta > + \int(e^{i<u,x>} - 1
- i
     < u,x > 1_{\{|x|\leq1\}})d\nu(x)]\}
\]
o\`u $\Sigma$ est une matrice $m\times m$ sym\'etrique semi-d\'efinie positive,
$\eta\in \RR^m$, et $\nu$ est une mesure positive sur $\RR^m$ v\'erifiant
$\nu\{0\} = 0$ et $\int 1 \wedgeÊ|x|^2 d\nu(x) < +\infty$. ($\nu$ est donc 
$\sig$-finie sur $\RR^m \backslash \{0\}$).

Dans cette formule, le triplet $(\Sigma,\eta,\nu)$ est d\'etermin\'e de fa\c{c}on
unique et s'appelle les caract\'eristiques locales de $Y$. Notons, cependant,
que d'autres formes peuvent \^etre donn\'ees \`a la formule de
L\'evy-Khintchine (cf [15]) et que le coefficient $\eta$ s'en trouve modifi\'e.
En revanche, $\Sigma$ et $\nu$ ont une signification intrins\`eque~: $Y$ est une
semi-martingale pour $(\Fc_\sig)$ et sa partie martingale continue $Y^c$ (nulle
en z\'ero) \`a pour crochet matriciel~:

\begin{equation}
\label{equa.5}
     < Y^c,Y^c >_\tau = \tauÊ\Sigma,
\end{equation}
et la mesure de L\'evy $\nu$ est d\'etermin\'ee par les sauts de $Y$~: la
mesure $d\tau \times d\nu$ sur $\RR_+ \times \RR^m$ est la projection
pr\'evisible duale de la mesure al\'eatoire~:

\[
     \mu(d\tau,dx) = \sum_{\sig>0} 1_{\{Y_\sig\neq Y_{\sig-}\}}
     \del_{(\sig,\Del Y_\sig)}(d\tau,dx).
\]
o\`u $\Del Y_\sig = Y_\sig-Y_{\sig-}$. C'est-\`a-dire qu'on a la formule:

\begin{equation}
\label{equa.6}
     \EE \sum_{\sig>0} H(\ome,\sig,\Del Y_\sigma)= \EE \int_0^\infty d\tau \int_{\RR^m} H(\ome,\tau,y) d\nu(y)
\end{equation}
pour tout processus $H$ pr\'evisible positif.

Un cas particulier important est celui des P.A.I.S. \`a valeurs r\'eelles
positives, ou subordinateurs. Ils sont caract\'eris\'es par (cf par exemple
[5]).

\[
     \forall \lam \geq 0, \tau \geq \sig \geq 0Ê\qquad \qquadÊ\qquad \qquad
\]
\begin{equation}
\label{equa.7}
     \EE[e^{-\lam(Y_\tau-Y_\sig)}|\Fc_\sig] = e^{-(\tau-\sig)\varp(\lam)}
\end{equation}

\[
     \varp(\lam) = b \lam + \int(1 - e^{-\lam x}) d\nu(x) \qquad \qquad
\]
$b \in \RR_+$, $\nu$ mesure positive sur $\RR_+$ telle que $\int(1 \wedge x)
d\nu(x) < +\infty$.

Les fonctions $\varp$ de la forme indiqu\'ee en (\ref{equa.7}) sont les
fonctions de Bernstein nulles en z\'ero, elles sont en correspondance
bi-univoque avec les semi-groupes de convolution de probabilit\'es sur $\RR_+$,
continus pour la convergence \'etroite (cf [1]).

\section{Visco\'{e}lasticit\'{e} et  P.A.I.S.}

\subsection{Cas d'un param\`etre observable}

Un mat\'eriau visco\'elastique lin\'eaire sans
vieillissement ne pr\'esentant qu'un param\`etre observable $q$
 a pour r\'eponse impulsionnelle (formule (6))~:

\begin{equation}
\label{equa.12}
     f(t) = L + Kt + \int (1 - e^{-tx})d\nu(x)
\end{equation}
o\`u $L$ et $K$  sont positifs et $\nu$ est une mesure de L\'evy de
subordinateur.

On d\'{e}finit une bijection si on associe \`{a} ce mat\'{e}riau le subordinateur $X$ tel que $X_0=L$
et dont la fonction de Bernstein est $f(t)-L$ (on pourrait aussi consid\'{e}rer un subordinateur \`{a} dur\'{e}e de vie
exponentielle pour traiter la constante $L$).

Remarquons que si on se donne un P.A.I.S. $Y$ r\'{e}el quelconque, $Y$ est une semi-martingale
 pour sa filtration naturelle. Soit $Y^c$ sa partie martingale continue. On peut associer
 \`{a} $Y$ un mat\'{e}riau visco-\'{e}lastique en posant
$Y_0^2=L$ et

\begin{equation}
\label{equa.13}
     f(t)=Y_0^2+\EE[tY_1^{c2}+ \sum_{0<\sig\leq 1}(1-e^{-t\Delta Y_\sigma^2})]
\end{equation}

\subsection{Cas de plusieurs param\`etres observables}

Consid\'erons un P.A.I.S. $Y$ \`a valeur $\RR^m$ tel que $Y_0$ soit une constante de $\RR^m$.
 On note $(\Sig,\eta,m)$ les caract\'eristiques de $Y$ (formule
(8)). $Y$ est une semi-martingale pour sa filtration naturelle
, soit $Y^c$ sa partie martingale continue.

On peut associer alors \`{a} $Y$ le mat\'{e}riau visco\'{e}lastique \`{a} $m$ param\`{e}tres observables dont les r\'{e}ponses impulsionnelles sont

\begin{equation}
\label{equa.14}
     f_{ij}(t)=Y_0^iY_0^j+\EE[tY_1^{ci}Y_1^{cj}+
 \sum_{0<\sig\leq 1}(1-e^{-t|\Delta Y_\sigma|^2})
\frac{\Delta Y_\sigma^i\Delta Y_\sigma^j}{|\Delta Y_\sigma|^2}]
\end{equation}
En effet cela donne

\begin{equation}
     f_{ij}(t)=Y_0^iY_0^j+t\Sigma_{ij}+\int_{\RR^m\backslash\{0\}}
(1-e^{-t|y|^2})\frac{y^iy^j}{|y|^2}\;dm(y)
\end{equation}
et les matrices $(\nu_{ij})$ o\`u $\nu_{ij}$ est
l'image par $y \rightarrow |y|^2$ de la mesure $\frac{y_i
y_j}{|y|^2}\cdot m(dy)$ sont toutes les matrices sym\'etriques
semi-d\'efinies positives de mesures sur $\RR^\ast_+$ v\'erifiant $\int_{\RR^\ast_+}
(x \wedge 1) d|\nu_{ij}|(x) < +\infty$.

\noindent{\bf Remarque} Notons que, si le mat\'eriau 1 [resp. 2] 
est associ\'{e} au processus $Y^{(1)}$
[resp. $Y^{(2)}$] et si $Y^{(1)}$ et $Y^{(2)}$ sont
pris ind\'ependants, le mat\'eriau dont le mod\`ele rh\'eologique est la mise
en s\'erie des mod\`eles des mat\'eriaux 1 et 2 (c'est-\`a-dire qui a les
m\^emes param\`etres observables et tel que

$
     f_{ij}(t) = f_{ij}^{1}(t) + f_{ij}^{2}(t))
$
est associ\'e au processus $Y^{(1)}+Y^{(2)}$.

\section{Mat\'{e}riaux conjugu\'{e}s et  dictionnaires}	

\subsection{Etude des fonctions de relaxation}

Pour voir \`a quoi correspond la mise en parall\`ele des mod\`eles
rh\'eologiques, consid\'erons les fonctions de relaxation, en nous pla\c{c}ant
d'abord dans le cas d'un seul param\`etre observable~:

Consid\'er\'ees comme distributions sur $\RR$ nulles pour $t < 0$, la r\'eponse
impulsionnelle $f$ et la fonction de relaxation sont li\'ees par la relation~:

\begin{equation}
\label{equa.20}
     f' * r' = \del_0.
\end{equation}
o\`{u} $*$ d\'{e}signe le produit de convolution.
 Notant $\hat{}$ la transformation de Laplace, la formule (6) donne~:

\begin{equation}
\label{equa.21}
     \widehat{(f')}(\theta) = L + \frac{K}{\theta} + \int_{\RR_+} \frac{x}{\theta
+ x}
     d\nu(x)Ê\qquad \theta > 0.
\end{equation}
Ceci nous conduit \`a la d\'efinition suivante~:

\begin{dfn}
Une fonction $h$ de $\RR^\ast_+$ dans $\RR$ est appel\'ee transform\'ee de
Stieltjes s'il existe $a \geq 0$ et une mesure $\mu$ positive sur $\RR_+$ tels
que~:

\[
     \forall \theta > 0 \qquad h(\theta) = a + \int \frac{1}{\theta + x}
     d\mu(x).
\]
\end{dfn}
Si $h$ est une transform\'ee de Stieltjes le couple $(a,\mu)$ est uniquement
d\'etermin\'e (cf [18]). La relation (\ref{equa.21}) exprime que $(f')\hat{}$
est transform\'ee de Stieltjes du couple ($L$, $K\del_0 + x\cdot\nu(x)$).

D'apr\`es le th\'eor\`eme 2 de [18], si $h$ est une transform\'ee de
Stieltjes non nulle, il en est de m\^eme de $\frac{1}{h(1/\theta)}$. Il
existe donc un couple $(\alp,\mu)$ unique non nul tel que~:

\[
     \frac{1}{\widehat{(f')}(\frac{1}{\theta})} = \alp + \int_{\RR_+}
     \frac{1}{\theta + x} d\mu(x)
\]
avec $\alp \geq 0$ et $\int_{\RR_+} \frac{1}{1 + x} d\mu(x) < +\infty$.

Par (\ref{equa.20}), on a donc~:

\[
     \widehat{(r')}(\theta) = \frac{1}{\widehat{(f')}(\theta)} = \alp + \int_\RR
     \frac{\theta}{1 + \theta x} d\mu(x)
\]
d'o\`u l'on tire~:

\[
     r(t) = \left(\alp + \mu\{0\}Ê\del_0 + \int_{]0}^{\infty}
\frac{e^{-\frac{t}{x}}}{x}
      d\mu(x)\right) 1_{\{t\geq0\}}.
\]
Ce qui \'etablit la :

\begin{pro}
Les fonctions de relaxation sont de la forme
\[
     r(t) = \alp + \bet \del_0 + \int_{\RR^\ast_+} e^{-tx} d\rho(x)
\]
o\`u $\alp, \bet \geq 0$ et $\rho$ est une mesure positive sur $\RR^\ast_+$ telle
que $\int_{\RR^\ast_+} \frac{1}{1 + x} d\rho(x) < +\infty$.
\end{pro}

En effet, il suffit de poser $\rho(dy) = y \tilde{\mu}(dy)$ o\`u $\tilde{\mu}$
est l'image de $\mu|_{\RR^\ast_+}$ par $x \rightarrow \frac{1}{x}$. La condition
$\int \frac{1}{1 + x} d\mu(x) <  +\infty$ impose alors~:

\begin{equation}
\label{equa.22}
     \int_{\RR^\ast_+} \frac{1}{1 + y} d\rho(y) < +\infty.
\end{equation}

\noindent{\bf Remarque.}
Ainsi les restrictions $t > 0$ des fonctions de relaxation ne sont pas toutes
les fonctions compl\`etement monotones, \`a cause de la condition
(\ref{equa.22}) (qui garantit que $r$ est une distribution au sens de Schwartz).

\subsection{Relation de conjugaison}

De m\^eme, dans le cas de $m$ param\`etres observables, les r\'eponses
impulsionnelles et les fonctions de relaxation sont reli\'ees par~:

\begin{equation}
\label{equa.23}
     (\sum_k f'_{ik} * r'_{k_j})_{ij} = \del_0\cdot I
\end{equation}
en notant $I$ la matrice identit\'e $m \times m$. Et par le m\^eme
raisonnement, on obtient que les fonctions de relaxation sont les distributions
de la forme

\begin{equation}
\label{equa.24}
     r_{ij}(t) = A_{ij} + B_{ij} \del_0 + \int_{\RR^\ast_+} e^{-tx} d\rho_{ij}(x)
\end{equation}
o\`u $(A_{ij}$, $(B_{ij})$ et $(\rho_{ij})$ sont des matrices sym\'etriques
semi-d\'efinies positives, les $\rho_{ij}$ \'etant des mesures telles que~:

\begin{equation}
\label{equa.25}
     \int_{]0}^{\infty} \frac{1}{1 + x} d|\rho_{ij}|(x) < +\infty.
\end{equation}
Il r\'esulte de (\ref{equa.24}) que {\it les primitives $R_{ij}(t)$ nulles pour $t <
0$ des $r_{ij}(t)$ sont exactement les fonctions de la forme (7). Et
par cons\'equent sont les r\'eponses impulsionnelles d'un mat\'eriau
visco\'elastique lin\'eaire sans vieillissement que nous appellerons
{\bf mat\'eriau conjugu\'e} du mat\'eriau initial}.

La relation de conjugaison s'exprime sur les r\'eponses impulsionnelles par

\begin{equation}
\label{equa.26}
     \sum_k f_{ik}^{(1)} * f_{k_j}^{(2)}(t) = \del_{ij}\cdot\frac{t^2}{2}
     1_{\{t\geq 0\}}.
\end{equation}

La relation de conjugaison est sym\'etrique et la mise en s\'erie des mod\`eles
rh\'eologiques de deux mat\'eriaux correspond \`a la mise en parall\`ele des
mat\'eriaux conjugu\'es et r\'eciproquement.\\

\noindent{\bf Dictionnaires et exemples}

Nous nous limitons ici, pour simplifier, au cas d'un seul param\`etre observable\\

{\small
\begin{tabular}{|l|c|l|}
     \hline
     & &\\
		   \qquad Mat\'eriaux & R\'eponses impulsionnelles & \qquad \qquad \quad
     processus\\
     & &\\
	    \hline\hline
     & &\\
		  mise en s\'{e}rie & addition des r\'eponses  & addition de
     processus\\
     & impulsionnelles & ind\'{e}pendants\\
  & &\\
	   \hline
     & &\\
		  mise en parall\`{e}le &addition des fonctions & addition des
     Processus\\
     & de relaxation & conjugu\'{e}s ind\'{e}pendants\\
  & &\\
	    \hline \hline
     & &\\
			  mat\'eriau & $f(t) = \frac{1}{a} 1_{\{t\geq0\}}$ & processus
           constant\\
     \'elastique & & de valeur $\frac{1}{a}$\\
	    & &\\
     \hline
     & &\\
     amortisseur & & translation uniforme \\
     (conjugu\'e du & $f(t) = at 1_{\{t\geq0\}}$ & de vitesse a\\
     pr\'ec\'edent) & &\\
     & &\\
     \hline \hline
     & &\\
     mat\'eriau de & $f(t) = (at + b) 1_{\{t\geq0\}}$ & translation uniforme
     de\\
     Maxwell & & vitesse $a$ issue du point $b$\\
     & &\\
     \hline
     & &\\
     mat\'eriau de & & processus de Poisson de saut\\
     Kelvin-Voigt & $f(t)= \frac{1}{a}(1-e^{-\frac{a}{b}t}) 1_{\{t\geq0\}}$ &
     d'amplitude $\frac{a}{b}$ et\\
     (conjugu\'e du & & d'intensit\'e $\frac{1}{a}$\\
     pr\'ec\'edent) & &\\
     & &\\
     \hline
     \end{tabular}
     
     \begin{tabular}{|l|c|l|}
      \hline
       & &\\
		   \qquad Mat\'eriaux & R\'eponses impulsionnelles & \qquad \qquad \quad
     processus\\
     \hline
     & &\\
     combinaison finie & & somme de processus de Poissons\\
     d'amortisseurs & & ind\'ependants de sauts \\
					et de ressorts & $f(t) = (\sum_i \alp_i(1-e^{-\lam_it})+at+b)
     1_{\{t\geq0\}}$ &  d'amplitude $\lam_i$ et d'intensit\'e $\alp_i$\\
	    & & et d'une translation de vitesse \\
     & & $a$ issue de $b$\\
     & &\\
     \hline
     & &\\
     combinaison finie & $f(t)=(\sum_i \bet_i(1-e^{-\mu_i t})+ct+d) 1_{t \geq 0}$ &
     somme de processus de Poissons\\
     d'amortisseurs et & & ind\'ependants de sauts \\
     de ressorts & $(\lam_1<\mu_1<\lam_2<\mu_2<\cdots$ si $c \neq 0$ &
     d'amplitude $\mu_i$ et d'intensit\'e $\bet_i$\\
     (conjugu\'e du & ou $\mu_1<\lam_1<\mu_2<\lam_2<\cdots$ si $a \neq 0$) & et
     d'un translation de\\
     pr\'ec\'edent) & & vitesse $c$ issue de $d$\\
     & &\\
     \hline
\end{tabular}}\\
\vspace{.5cm}

\noindent La relation de conjugaison attire l'attention sur le cas suivant\\

\begin{center}
\begin{tabular}{|l|l|}
      \hline
      & \\
      \qquad R\'eponses impulsionnelles & \qquad \qquad  Processus\\
      & \\
      \hline \hline
      &Ê\\
      $f(t) = \frac{ t^\alp }{\alp\Gam(\alp)}1_{\{t\geq0\}} \quad \alp \in
      ]0,1[$ & processus stable unilat\'eral d'ordre $\alp$\\
      & \\
      conjugu\'e avec le mat\'eriau & \\
      d'indice $1-\alp$ &Ê\\
      &\\
      \hline \hline
      & \\
	     seul mat\'eriau qui est son  propre conjugu\'e & processus stable unilat\'eral d'ordre $\Frac$\\
      & \\
      $f(t) = 2 \sqrt{\frac{t}{\pi}} 1_{\{t\geq0\}}$ & \\
      & \\
      \hline
\end{tabular}\\
\end{center}

\hspace{.5cm}

\noindent{\bf Autres remarques.}\\

1 - L'importante propri\'et\'e du c\^one des fonctions de Bernstein d'\^etre
stable par composition (cf [1]) fait que si $f_1(t)$ et $f_2(t)$ sont prises
dans les tableaux ci-dessus, $f_1((f_2(t))$ est encore la r\'eponse
impulsionnelle d'un mat\'eriau visco\'elastique lin\'eaire sans vieillissement,
ce qui fournit une grande vari\'et\'e de fonctions.

2 - La r\'eponse impulsionnelle d'un mat\'eriau visco\'elastique \`a un 
param\`etre observable peut \^etre consid\'er\'ee comme une horloge. C'est, en
effet, un ph\'enom\`ene d\'eclench\'e \`a l'instant $0$ et qui \'evolue ensuite
ind\'efiniment selon une dynamique propre. Ces horloges ont la propri\'et\'e
suivante~:

Consid\'erons :

a) un mobile $A$ anim\'e d'une vitesse rectiligne uniforme, en z\'ero \`a
l'instant $0$.

b) deux curseurs $B$ et $C$  li\'es respectivement aux r\'eponses impulsionnelles
de mat\'eriaux $(1)$ et $(2)$ et pouvant se d\'eplacer parall\`element \`a $A$.

Si on note les positions successives $a_0 = 0,a_1,\cdots,a_n,\cdots$ de $A$ 
lorsque $B$ atteint des distances r\'eguli\`eres $0,h,2h,\cdots,nh,\cdots,$ 
alors les positions successives $\alp_0 = 0, \alp_1,\cdots,\alp_n,\cdots$ de
$A$ lorsque $C$ atteint les $a_n$ sont aussi celles qu'il atteindrait
lorsqu'un corps visco\'elastique $D$ voit son index atteindre les distances
$nh$. Cette propri\'et\'e est due \`a la stabilit\'e par composition du c\^one
des fonctions de Bernstein.

3 - La th\'eorie de Biot a fait l'objet de nombreuses extensions, d'abord au
cas avec vieillissement, \'egalement avec des concepts de g\'eom\'etrie
diff\'erentielle les param\`etres cach\'es \'etant suppos\'es tensoriels cf
[22] et \`a des cas non lin\'eaires cf [17]. Il est vraisemblable qu'une partie
de ces extensions peut recevoir une interpr\'etation \`a partir des
caract\'eristiques locales des semi-martingales.\\

\noindent{\bf Bibliographie}\\

\noindent[1] Ch. Berg, G. Forst
      {\it Potential theory on locally compact abelian groups}.\\
	    Springer  (1975).

\noindent[2] M. Biot
     Theory of stress-strain relations in anisotropic viscoelasticity and
     relaxation phenomena.\\
		   J. Appl. Phys. 25, 11, 1385-1391 (1954).

\noindent[3] M. Biot
      Thermoelasticity and irreversible thermodynamics.
     {\it J. Appl. Phys.}, 27, 3, 240-253 (1956).

\noindent[4] J.M. Bony, Ph. Courr\`{e}ge, P. Priouret,
     Semigroupes de Feller sur une vari\'{e}t\'{e} \`{a} bord compacte
     {\it Ann. Inst. Fourier}, 18, 369-521, (1968)

\noindent[5] N. Bouleau
      {\it Processus stochastiques et applications}.
	    Hermann (1988).

\noindent[6] N. Bouleau
					Interpr\'{e}tation probabiliste de la visco\'{e}lasticit\'{e}
 lin\'{e}aire,  {\it mechanics Research Comm.} vol 19, 16-20, (1992)

\noindent[7] N. Bouleau et O. Chateau, 
     Le processus de la subordination, {\it Note C. R. Acad. Sci. Paris}. t. 309, sI, 955-958, (1989)

\noindent[8] O. Chateau, 
"Quelques remarques sur les PAIS et la subordination au sens de Bochner",
 Th\'{e}se Univ. Paris VI, (1990)

\noindent[9] B.D. Coleman
      On thermodynamics, strain impulses and viscoelasticity.
	     {\it Arch. Rat. Mech. Anal.} Vol. 17, 230-254 (1964).

\noindent[10] C. Dellacherie, P.A. Meyer
      {\it Probabilit\'es et Potentiel},
	    Chapitres V \`a VIII, Th\'eorie des martingales, Hermann(1980).

\noindent[11] C. Dellacherie, P.A. Meyer 
     {\it Probabilit\'es et potentiel}, chapitres XII \`a XVI,
     Hermann (1987).

\noindent[12] R. Durrett
	    {\it Brownian motion and martingales in analysis},
     Wadsworth (1984).

\noindent[13] C. Eckart
     Theory of elasticity and anelasticity,
     {\it Phys. Rev.} Vol. 73, 373-382 (1948).

\noindent[14] S.M. Ermakov, V.V. NekrutkinN, A.S. Sipin
     {\it Random processes for classical equations of math\'ematical physics},\\
     Kluwer (1989).

\noindent[15] W. Feller
      {\it Introduction to probability theory and its applications}.
     Tome II, Wiley (1966).

\noindent[16] Y.C. Fung,
      {\it Foundations of solid mechanics}.
      Prentice Hall (1965).

\noindent[17] B. Halphen, N. Guyen Quoc-Son
     Sur les mat\'eriaux standard g\'en\'eralis\'es.
	    {\it J. M\'ecanique}, 14, 1 39-63 (1975).

\noindent[18] F. Hirsch,
     Int\'egrales de r\'esolvantes et calcul symbolique.
	    {\it Ann. Inst. Fourier} XXII, 4, 239-264 (1972).

\noindent[19] J. Jacod,
      {\it Calcul stochastique et probl\`emes de martingales}.
     Lecture Notes in Math. 714, Springer (1979).

\noindent[20] J. Mandel,
      {\it M\'ecanique des milieux continus}.
      Tome II, Gauthier-Villars (1966).

\noindent[21] J. Meixmer,
     Die thermodynamische theorie der Relaxationser-Scheinunger und ihr
     zuzammenhang mit der Nachwirkungstheorie,
	    {\it Kolloid Z.} Vol. 14, 2-9 (1953).

\noindent[22] F. Sidoroff,
      {\it Journal de M\'ecanique},
      Vol. 14, 3 (1975).

\noindent[23] K. Yosida,
				 {\it Functional Analysis }. Springer 1974.

\end{document}